\documentclass[12pt]{article}

\usepackage{latexsym}
\usepackage{amsmath}
\usepackage{amsfonts}
\usepackage{amssymb}

\newcommand{\bean}{\begin{eqnarray}}
\newcommand{\eean}{\end{eqnarray}}
\newcommand{\bea}{\begin{eqnarray*}}
\newcommand{\eea}{\end{eqnarray*}}
\newcommand{\bsa}{\begin{subarray}{c}}
\newcommand{\esa}{\end{subarray}}
\newcommand{\bi}{\begin{itemize}}
\newcommand{\ei}{\end{itemize}}

\newtheorem{lemma}{Lemma}[section]
\newtheorem{thm}[lemma]{Theorem}

\begin{document}

\title{\bf Structure of the module of vector-valued modular forms}
\author{Christopher Marks\thanks{Supported by NSA} \\
Geoffrey Mason\thanks{Supported by NSA and NSF}\\
University of California at Santa Cruz} 
\date{}
\maketitle

\abstract{Let $V$ be a representation of the modular group $\Gamma$ of dimension $p$. We show that
the $\mathbb{Z}$-graded space $\mathcal{H}(V)$ of holomorphic vector-valued modular forms
associated to $V$ is a free module of rank $p$ over the algebra $\mathcal{M}$ of classical holomorphic modular forms. We study the nature of $\mathcal{H}$ considered as a functor from $\Gamma$-modules to graded $\mathcal{M}$-lattices and give some applications, including the calculation of the 
Hilbert-Poincar\'{e} of $\mathcal{H}(V)$ in some cases.}\\
MSC 11F99, 13C05

\section{Introduction}
Vector-valued modular forms have been a part of number theory for some time,  but a systematic development of their properties has begun only relatively recently  (\cite{BG1}, \cite{KM1}-\cite{KM3},
\cite{M1}, \cite{M2}). One motivation for this
comes from rational and logarithmic field theories, where vector-valued modular forms arise naturally(\cite{ES}, \cite{DLM}, \cite{My}, \cite{Z}). Modular forms on noncongruence subgroups are a special case of vector-valued modular forms, and one of the goals of both this case and the general
theory  is to find arithmetic conditions which characterize \emph{classical} modular forms
(i.e. on a congruence subgroup) among all vector-valued modular forms (cf. \cite{KoM}, \cite{KL})
A systematic approach to this problem requires detailed information about the general structure of the space of vector-valued modular forms attached to a representation of the modular group. The purpose of this paper is to provide such information. The main theorems extend some of the results of \cite{M2}, dealing with $2$-dimensional irreducible representations, to a general context.
In the next few paragraphs we give some
 basic definitions sufficient to state our main results. A fuller discussion of background material can be found in \cite{M2} and in Section 2 below.

\bigskip
Let $\Gamma = SL(2, \mathbb{Z})$,  acting on the complex upper half-plane $\mathcal{H}$ in the usual way, and with standard generators
\begin{eqnarray*}
S = \left(\begin{array}{cc}0 & -1 \\ 1 & 0\end{array}\right) , \ \ T = \left(\begin{array}{cc} 1 & 1 \\0 & 1\end{array}\right).
\end{eqnarray*}
Let $V$ be a left $\mathbb{C} \Gamma$-module of finite dimension $p$. We denote the action of
$\gamma \in \Gamma$ on $v \in V$ by $\gamma v$ or $\gamma. v$.  Let
\begin{eqnarray*}
\frak{F} = \{\mbox{holomorphic} \ F: \mathcal{H} \rightarrow V \}. 
\end{eqnarray*}
For an integer $k$ there is a right action $\frak{F} \times \Gamma \rightarrow \frak{F}$ given by
\begin{eqnarray}\label{FGact}
\gamma: F(\tau)  \mapsto \gamma^{-1}.F|_k \gamma (\tau),
\end{eqnarray}
 where $|_k$ is the usual stroke operator familiar from the classical theory of modular forms. 
 A \emph{weak}  vector-valued modular form
of weight $k$ is a $\Gamma$-invariant of the action (\ref{FGact}).

\medskip
There is a basic subdivision of the general theory according to whether or not the restriction of 
$V$ to the $T$-matrix
is unitary. If this is the case, we say that $V$ is  \emph{$T$-unitarizable}. We will always assume throughout the present paper that $V$ is indeed $T$-unitarizable. Then $V$ has a basis which
furnishes a representation $\rho: \Gamma \rightarrow GL(p, \mathbb{C})$ with 
\begin{eqnarray}\label{Tdiag}
\rho(T) = diag(e^{2 \pi i m_1}, \hdots, e^{2 \pi i m_p})
\end{eqnarray}
and $0 \leq m_j < 1$ for $1 \leq j \leq p$.
A weak modular form is then a column vector  of holomorphic \emph{component functions}\footnote{superscript $t$ denotes \emph{transpose}
of a vector} 
$F(\tau) = (f_1(\tau), \hdots, f_p(\tau))^t$ satisfying
\begin{eqnarray}\label{newdef}
\rho(\gamma)F(\tau) = (f_1|_k \gamma (\tau), \hdots, f_p|_k \gamma (\tau))^t, \ \ \gamma \in \Gamma.
\end{eqnarray}
By a standard argument,   (\ref{newdef}) implies that each component function has a $q$-expansion
\begin{eqnarray}\label{compqexp}
f_j(\tau) = \sum_{n \in \mathbb{Z}} a_{jn}q^{\lambda_j+n}
\end{eqnarray}
such that $\lambda_j - m_j$ is a nonnegative integer. We call $F(\tau)$ a (holomorphic) \emph{vector-valued modular form of weight $k$}
(with respect to $V$, or $\rho$) in case $a_{jn} = 0$ for  $n<0$ and $1 \leq j \leq p$. This condition is independent of the choice of basis of $V$.

\medskip
Let $\mathcal{H}(k, V)$ be the set of all vector-valued modular forms of weight $k$ with respect to $V$.
It is a finite-dimensional linear space which reduces to $0$ for $k \ll 0$. 
(See \cite{KM1}, \cite{KM2}, \cite{M1} for further discussion.)
In the further development of the theory it is no loss to assume that $\rho(S^2) = \pm I_p$
(cf. \cite{KM1}). With this assumption, the direct sum of the spaces $\mathcal{H}(k, V)$ takes the form
\begin{eqnarray}\label{graded}
\mathcal{H}(V) = \bigoplus_{k \geq 0} \mathcal{H}(k_0+2k, V)
\end{eqnarray}
for a certain minimal weight $k_0$. 

\medskip
Let $\mathcal{M} = \oplus_k \mathcal{M}_{2k}$ be the space of (classical holomorphic) modular forms on $\Gamma$, 
regarded as a $\mathbb{Z}$-graded, weighted polynomial algebra.
Pointwise multiplication by elements of $\mathcal{M}$ turns $\mathcal{H}(V)$ into a $\mathbb{Z}$-graded $\mathcal{M}$-module. 
We can now state our main results, which concern the structure of $\mathcal{H}(V)$
considered as $\mathcal{M}$-module. We emphasize that throughout the paper we assume that $V$ is $T$-unitarizable
of dimension $p$ with $\rho(S^2)= \pm I_p$.
 The basic result is

\medskip
\noindent
{\bf Theorem 1}\label{thm1} $\mathcal{H}(V)$ is a free
$\mathcal{M}$-module of rank $p$.

\bigskip
We may therefore consider $\mathcal{H}$ as a \emph{covariant functor}
\begin{eqnarray*}
\mathcal{H}: \Gamma\mbox{-}\mbox{Mod}_{fin} \rightarrow \mbox{Gr}\mathcal{M}\mbox{-Latt}
\end{eqnarray*}
from the category of finite-dimensional $\mathbb{C}\Gamma$-modules to the category
of $\mathbb{Z}$-graded $\mathcal{M}$-lattices (an $\mathcal{M}$-lattice is a finitely generated
projective ($=$ free) $\mathcal{M}$-module). The next Theorem gives some of the main properties of 
this functor.

\medskip
\noindent
{\bf Theorem 2}\label{thm2} Suppose that 
$0 \rightarrow U \stackrel{i}{\rightarrow} V \stackrel{j}{\rightarrow} W \rightarrow 0$
is a short exact sequence of $\mathbb{C}\Gamma$-modules. Then the following hold.\\
(a) $\mathcal{H}$ is faithful and left exact, in particular there is an exact sequence
\begin{eqnarray*}
0 \rightarrow \mathcal{H}(U) \stackrel{i_*}{\rightarrow} \mathcal{H}(V) \stackrel{j_*}{\rightarrow} 
\mathcal{H}(W).
\end{eqnarray*}
(b) $\mathcal{H}$ is \emph{pure} in the sense that $i_* \mathcal{H}(U)$ is an 
$\mathcal{M}$-direct summand of $\mathcal{H}(V)$.\\
(c) $j_*\mathcal{H}(V)$ is a free $\mathcal{M}$-module of rank $\dim W = \dim V$-$\dim U$
and it has \emph{finite codimension} in $\mathcal{H}(W)$.

\bigskip
Let $D: \mathcal{M}_{k} \rightarrow \mathcal{M}_{k+2}$
be the usual graded derivation of $\mathcal{M}$ defined via
\begin{eqnarray}\label{Dact}
D(f) = \theta f + kPf, \ f \in \mathcal{M}_{2k}, 
\end{eqnarray}
where $\theta = qd/dq = (2\pi i)^{-1}d/d\tau$ and
\begin{eqnarray*}
P = E_2(\tau) = -1/12 + \sum_{n \geq 1} \sigma_1(n)q^n.
\end{eqnarray*}
The algebra $\mathcal{R} = \mathcal{M}[d]$ of skew polynomials in $d$ with coefficients in 
$\mathcal{M}$ consists of (noncommutative) polynomials
$\sum_{j = 0}^n f_jd^j,  f_j \in \mathcal{M}$, 
with the usual addition and multiplication subject to the identity
\begin{eqnarray*}
df - fd = D(f).
\end{eqnarray*}
The action of $\mathcal{M}$ on $\mathcal{H}(V)$ extends to an action of $\mathcal{R}$, so
that $\mathcal{H}(V)$ is a graded $\mathcal{R}$-module. See \cite{M1}, \cite{M2} for further details. Although this $\mathcal{R}$-module
structure is not explicit in the statements of Theorems 1 and 2, it plays an important r\^{o}le
in the proofs. 

\bigskip
One of the main steps in the proof of Theorem \ref{thm1} is to establish that
$\mathcal{H}(V)$ is \emph{finitely generated} as $\mathcal{M}$-module. The theory of vector-valued Poincar\'{e} series \cite{KM1} implies the existence $F \in \mathcal{H}(k, V)$ such that the component functions of $F$ are linearly independent, and together with the theory of differential equations implies
that $\mathcal{R}F$ contains a free graded $\mathcal{M}$-submodule $\mathcal{H}'$ of rank $p$
such that $\mathcal{H}(V)/\mathcal{H}'$ has bounded degree. Methods of commutative algebra then imply  finite generation. Freeness is established by showing that
$\mathcal{H}(V)$ is Cohen-Macaulay, which turns out to be very natural in the present context.

\bigskip
Theorem 1 implies that there are 
$p$ distinguished integers $e_1, \hdots, e_p$, namely the weights of a set
of $p$ vector-valued modular forms $F_1, \hdots, F_p$ which are free generators
of $\mathcal{H}(V)$ as $\mathcal{M}$-module, that  are \emph{uniquely determined}
by $V$. We call these the \emph{fundamental weights} associated to $V$. 
In terms of  the Hilbert-Poincar\'{e} series
\begin{eqnarray*}
PS \ \mathcal{H}(V) = \sum_{k\geq0} \dim \mathcal{H}(k_0+2k)t^{k_0+2k}
\end{eqnarray*}
we have
\begin{eqnarray}
PS \ \mathcal{H}(V) = \frac{t^{e_1}+ \hdots + t^{e_p}}{(1-t^4)(1-t^6)}.
\end{eqnarray}

\bigskip
Exactly how the fundamental weights are determined by $V$ remains somewhat mysterious.
To some extent, this is related to the fact that $\mathcal{H}$ is \emph{not} an exact functor
(cf. Theorem 4 below). Nevertheless, we can use the fundamental weights to impose a $\mathbb{Z}$-grading on $V$, say by choosing a basis
$v_i$ and giving $v_i$ weight $e_i$. Then the functor $\mathcal{H}$ corresponds to an extension of scalars
\begin{eqnarray*}
V \mapsto \mathcal{M} \otimes_{\mathbb{C}} V
\end{eqnarray*}
where $\mathcal{M} \otimes V$ inherits the tensor product grading.

\medskip
We prove two further results, based on Theorems 1 and 2, that illustrate aspects of
the general theory. The first is concerned with the case that $\mathcal{H}(V) = \mathcal{R}F$ is a 
\emph{cyclic} $\mathcal{R}$-module. This condition necessarily holds in a number of cases when $V$ is an \emph{irreducible}
$\mathbb{C}\Gamma$-module of small dimension, including all irreducible $V$ with $\dim V \leq 3$.
 (See \cite{M2} and \cite{Ma} for further details.)
We establish

\medskip
\noindent
{\bf Theorem 3}\label{thm3} 
If $\mathcal{H}(V) = \mathcal{R}F$ is a \emph{cyclic} $\mathcal{R}$-module with generator $F$
of weight $k_0$, then the component functions of $F$ form a fundamental system of solutions of a
\emph{modular linear differential equation} (MLDE)
\begin{eqnarray}\label{MLDE}
L_{k_0}[f] = 0
\end{eqnarray}
of weight $k_0$ and order $p$. The roots of the indicial equation
are the exponents $m_j$ in (\ref{Tdiag}), they are \emph{distinct}, and the minimal weight $k_0$ satisfies
\begin{eqnarray}\label{wtrestr}
12\sum_{j=1}^p m_j = p(p+k_0-1).
\end{eqnarray}
Conversely, if the indicial equation of the MLDE (\ref{MLDE}) has real, distinct roots
$m_j$  that lie between $0$ and $1$ and  satisfy (\ref{wtrestr}),  then a fundamental system of solutions
spans a $\Gamma$-module $V$ and are the components of a vector-valued modular form
$F$ for which $\mathcal{H}(V) = \mathcal{R}F$.
(For further information concerning MLDE's in this context, see \cite{M1}, \cite{M2}.)

\medskip
In this case we can take the free generators to be $F, DF, \hdots, D^{p-1}F$, whence the fundamental weights are $k_0, k_0+2, \hdots, k_0+2p-2$, and
\begin{eqnarray*}
PS \ \mathcal{H}(V) = \frac{t^{k_0}(1-t^{2p})}{(1-t^2)(1-t^4)(1-t^6)}.
\end{eqnarray*}

\medskip
We complete the paper by discussing the case of \emph{indecomposable}  $2$-dimensional representations of $\Gamma$. The irreducible case is handled
 in \cite{M2} (alternatively, by Theorem 3). The case when $V$ is indecomposable but
 \emph{not} irreducible is less straightforward, and  illustrates some of the 
subtleties involved in calculating the Hilbert-Poincar\'{e} series in general. 

\medskip
To state our result, recall
that the group of characters of $\Gamma$ is cyclic of order $12$,
generated by a character $\chi$ 
 uniquely determined by the equality $\chi(T) = e^{2\pi i/12}$.

\medskip
\noindent
{\bf Theorem 4}\label{thm2dindecomp} Suppose that $V$ is a $2$-dimensional indecomposable
$\Gamma$-module occurring in the short exact sequence
\begin{eqnarray*}
0 \rightarrow \chi^a \stackrel{i}{\rightarrow} V \stackrel{j}{\rightarrow} \chi^b \rightarrow 0
\end{eqnarray*}
($0 \leq a, b \leq 11, \ |a-b|=2$) furnishing a representation $\rho$ which is \emph{upper triangular}.
 One of the following holds:
 
 \medskip
 \noindent
(a) There is a (split) short exact sequence of graded $\mathcal{M}$-modules
\begin{eqnarray*}
0 \rightarrow \mathcal{H}(\chi^a) \stackrel{i_*}{\rightarrow} \mathcal{H}(V) \stackrel{j_*}{\rightarrow} \mathcal{H}(\chi^b) \rightarrow 0,
\end{eqnarray*}
(b) $(a, b) = (10, 0)$ or $(11, 1)$ and there is an exact sequence of graded $\mathcal{M}$-modules
\begin{eqnarray*}
0 \rightarrow \mathcal{H}(\chi^a) \stackrel{i_*}{\rightarrow} \mathcal{H}(V) \stackrel{j_*}{\rightarrow} \mathcal{H}(\chi^b) \rightarrow \mathbb{C}_{b} \rightarrow 0
\end{eqnarray*}
($\mathbb{C}_b$ is the $1$-dimensional graded $\mathcal{M}$-module in weight $b$).

\medskip
In particular, part (b) confirms our earlier assertion that $\mathcal{H}$ is generally not right exact.
Note that the condition $|a-b|=2$ in Theorem 4 \emph{necessarily} holds whenever
$V$ is indecomposable but not irreducible. See Section 3 of \cite{M2}
and Section 4 below for further details.

\medskip
Terry Gannon has recently informed the authors that he and Peter Bantay have also found a proof
of Theorem 1 \cite{BG2}. Their methods, as in \cite{BG1},  are rather different to ours, and will appear elsewhere.

\section{Preliminaries}
We keep the assumptions and notation of  Section 1. The algebra of classical modular forms on 
$\Gamma$ is a $\mathbb{Z}$-graded
algebra
\begin{eqnarray*}
\mathcal{M} = \bigoplus_{k = 0}^{\infty} \mathcal{M}_{2k}
\end{eqnarray*}
where
$\mathcal{M}_{2k}$ is the space of holomorphic modular forms of weight $2k$.
There is an isomorphism of $\mathbb{Z}$-graded algebras
\begin{eqnarray*}
\mathcal{M} \stackrel{\cong}{\rightarrow} \mathbb{C}[Q, R]
\end{eqnarray*}
where $\mathbb{C}[Q, R]$ is a weighted polynomial algebra with generators
\begin{eqnarray*}
Q &=& E_4(\tau) = 1 +240\sum_{n \geq 1} \sigma_3(n)q^n, \\
R &=& E_6(\tau) = 1 - 504\sum_{n \geq 1} \sigma_5(n)q^n, 
\end{eqnarray*}
 the usual Eisenstein series of weights $4$ and $6$ respectively.
The first two Lemmas below are consequences of the theory of vector-valued Poincar\'{e} series  
\cite{KM1} and play an important r\^{o}le in the proofs of the main Theorems.

\begin{lemma}\label{lemmamus} Suppose that $\mu_1, \hdots, \mu_p$ is a sequence of \emph{nonnegative}
integers and $(c_1, \hdots, c_p)$ a sequence of scalars. Then for all large enough $k$, there is 
$F(\tau) = (f_1(\tau), \hdots, f_p(\tau))^t \in \mathcal{H}(k_0+2k, V)$
such that
\begin{eqnarray*}
f_r(\tau) = c_rq^{\mu_r+m_r} + \hdots, \ 1 \leq r \leq p.
\end{eqnarray*}
\end{lemma}
\begin{pf}  Choose $r$ in the range $1 \leq r \leq p$ and an integer $\nu_r < 0$. By Theorem 3.2 
of \cite{KM1} we can find,  for large enough $k$, a (meromorphic) vector-valued modular form  $P_r(\tau)$ of weight  $k$ such that  the  component functions $P_{r, j}(\tau)$ are as follows:
\begin{eqnarray*}
P_{r, r}(\tau) &=& q^{\nu_r + m_r} + \hdots   \\
P_{r, j}(\tau) &=& q^{n_{r,j} + m_j} + \hdots, \ n_{r,j}>\mu_j, \ j \not= r.
\end{eqnarray*}
We may, and shall, arrange that the
$P_{r}(\tau)$ have a \emph{common weight} $k$ for large enough $k$. Then
$P(\tau) = \sum_r c_rP_{r}(\tau)$ is a meromorphic vector-valued modular form of weight $k$. 
Now choose $\nu_r = \mu_r - d$ for some integer $d$. The vector-valued modular form 
\begin{eqnarray*}
F(\tau) = \Delta^d(\tau)P(\tau).
\end{eqnarray*}
 has the required properties, and the Lemma follows. $\hfill \Box$
\end{pf}

\medskip
We say that the vector-valued modular form $F(\tau)$
is \emph{essential} if $F(\mathcal{H}) \subseteq V$ \emph{spans} $V$. Choose a basis of $V$ and write $F(\tau)$ in component form, say
\begin{eqnarray*}
F(\tau) = (f_1(\tau), \hdots, f_p(\tau))^t.
\end{eqnarray*}
Then $F(\tau)$ is essential if, and only if, $f_1(\tau), \hdots, f_p(\tau)$ are linearly independent functions.
If, in Lemma \ref{lemmamus}, we choose the $\mu_r$ to be \emph{distinct}, the resulting
vector-valued modular form $F(\tau)$ has component functions which are clearly linearly independent.
Hence, we obtain

\begin{lemma}\label{essvvform} $\mathcal{H}(V)$ contains an essential vector-valued modular form.
$\hfill \Box$
\end{lemma}

\begin{lemma}\label{lemmaHest} There is a constant $A$ depending only on $V$ such that 
\begin{eqnarray}
\left| \dim \mathcal{H}(k_0+2k, V) - \frac{pk}{6} \right| < A
\end{eqnarray}
for all $k \geq k_0$ with $k$-$k_0$ even.
\end{lemma}
\begin{pf} This too follows from the theory of Poincar\'{e} series. Indeed, it is a consequence
of Theorems 2.5 and 4.2 of \cite{KM1}. The result may also be deduced from Lemma \ref{essvvform} using the theory of MLDEs as in \cite{M1}. $\hfill \Box$
\end{pf}

\bigskip
 Suppose that $\alpha: V \rightarrow W$ is a morphism of $\Gamma$-modules.
 The commuting diagram
\begin{eqnarray*}
 &&V \ \ \  \stackrel{\alpha}{\longrightarrow} \ \ \  W  \\
 && \ \ \   \nwarrow \ \ \ \ \ \nearrow \\
 &&\ \ \ \ \ \  \ \  \frak{H}
  \end{eqnarray*}
  allows us to push-foward holomorphic maps $F: \frak{H} \rightarrow V$ to get holomorphic maps
  $ \alpha \circ F: \frak{H} \rightarrow W$.  
  \begin{lemma}\label{lemma6.1}  The following hold: \\
  (a) \  $\alpha$ induces a map
  \begin{eqnarray*}
&&\alpha_*: \mathcal{H}(k, V) \rightarrow \mathcal{H}(k, W) \\
&& \ \ \ \ \ \ \  \ \ \ \ \    F \mapsto \alpha \circ F.
\end{eqnarray*}
(b) \ If \ $V \stackrel{\alpha}{\rightarrow} W \stackrel{\beta}{\rightarrow} X$ are morphisms
of $\Gamma$-modules, then $(\beta \circ \alpha)_* = \beta_* \circ \alpha_*$.
  \end{lemma}
  \noindent
  \begin{pf} Suppose that $F(\tau) \in \mathcal{H}(k, V)$. Then
  \begin{eqnarray*}
&&\gamma^{-1}. \alpha_*(F) |_k \gamma (\tau) =  
\gamma^{-1}.( \alpha \circ F)|_k \gamma (\tau) \\
= && 
 \alpha \gamma^{-1}. (F|_k \gamma (\tau)) \\
=&&\alpha \circ F(\tau) = \alpha_*(F)(\tau).
\end{eqnarray*}
This shows that $\alpha_*(F) \in \mathcal{H}(k, W)$, and part (a) follows. Part (b) is clear. $\hfill \Box$
  \end{pf}

\bigskip
From Lemma \ref{lemma6.1} it follows that there is a covariant functor
\begin{eqnarray*}
\mathcal{H}: \Gamma\mbox{-}Mod_{fin} \rightarrow Gr\mathcal{M}\mbox{-Mod}
\end{eqnarray*}
from the category of finite-dimensional $\Gamma$-modules to the category of
$\mathbb{Z}$-graded $\mathcal{M}$-modules.

\begin{lemma}\label{lemmaexact} The functor $\mathcal{H}$ is \emph{faithful} and \emph{left exact}.
\end{lemma}
\begin{pf} To prove that $\mathcal{H}$ is faithful, we must show that if
$\alpha: V \rightarrow W$ is a nonzero morphism of $\Gamma$-modules then
$\mathcal{H}(\alpha) = \alpha_*$ is also nonzero. By Lemma \ref{essvvform} there is an integer
$k$ such that $\mathcal{H}(k, V)$ contains an essential vector-valued modular form, say $F(\tau)$.
Then $F(\tau)$ generates $V$ as linear space, so if $\alpha \not= 0$ then also $\alpha_*(F) = \alpha \circ F \not= 0$. Thus $\alpha_*$ is itself nonzero, as required.

\bigskip
As for left exactness, we have to show that if 
\begin{eqnarray*}
0 \rightarrow U \stackrel{\alpha}{\rightarrow} V \stackrel{\beta}{ \rightarrow} W
\end{eqnarray*}
is exact, then so too is
\begin{eqnarray}\label{exseq}
0 \rightarrow \mathcal{H}(U) \stackrel{\alpha_*}{ \rightarrow} \mathcal{H}(V)
 \stackrel{\beta_*}{\rightarrow} \mathcal{H}(W).
\end{eqnarray}
This is standard, and we omit the proof.  $\hfill \Box$
\end{pf}

\bigskip
By a \emph{$d$-ideal} in $\mathcal{M}$ we mean a left $\mathcal{R}$-submodule of
$\mathcal{M}$ regarded as $\mathcal{R}$-module. In other words, an ideal in $\mathcal{M}$
 invariant under $D$.

\begin{lemma}\label{lemma3.1} Let $I \subseteq \mathcal{M}$ be a nonzero, graded
$d$-ideal. Suppose that
\begin{eqnarray}\label{Dcondition}
\mbox{Ann}_{\mathcal{M}/I}(\Delta) = 0.
\end{eqnarray}
Then $I = \mathcal{M}$.
\end{lemma}
\noindent
\begin{pf} Since $I$ is graded it is the direct sum of its subspaces 
$I_{2k} = I \cap \mathcal{M}_{2k}$. Let $m \in I_{2k}$ be a \emph{nonzero}  element of least weight
in $I$, and consider the linear span
$N$ of the two elements $QD(m), Rm \in I_{2k+6}$. If $\dim N = 2$ 
 then it contains a nonzero cuspform $\alpha$. Then $\alpha = \Delta \beta \in I$ 
 with 
 $\beta \in \mathcal{M}_{2k-6}$, and by
(\ref{Dcondition}) we have $\beta \in I$. This contradicts the minimality of $2k$, and shows that
$QD(m)$ and $Rm$ are scalar multiples of one another.
Set
\begin{eqnarray*}
m = \sum c_{uv}Q^uR^v
\end{eqnarray*}
where $c_{uv}$ is a scalar and $(u, v)$ range over pairs satisfying $2k = 4u + 6v$.
Then we must have
\begin{eqnarray*}
QD(m) &=& \sum c_{uv} (uQ^{u}R^{v+1} + vQ^{u+3}R^{v-1}) \\
&=& c \sum c_{uv}Q^uR^{v+1}
\end{eqnarray*}
for a nonzero constant $c$. If $u_0$ is the highest  power of $Q$ occurring with nonzero 
coefficient in the expression for $m$, we see that the corresponding $v$ is zero, and we have
\begin{eqnarray*}
c = u_0.
\end{eqnarray*}
Then looking at the \emph{lowest} power of $Q$ that occurs with nonzero coefficient, 
say $u_1$, we obtain
\begin{eqnarray*}
cc_{u_1v_1} = u_0c_{u_1v_1} = c_{u_1v_1}u_1.
\end{eqnarray*}
We conclude that $u_0 = u_1$, i.e. $m = cQ^{c}$. Then $D(m) = c^2Q^{c-1}R$ and $I$
contains 
\begin{eqnarray*}
3RD(m) - 2cQ^2m = c^2 Q^{c-1}(3R^2 - 2Q^3) = c^2Q^{c-1}\Delta.
\end{eqnarray*}
If $c \geq 1$, (\ref{Dcondition}) tells us that $Q^{c-1} \in I_{2k-4}$, contradiction. So $c=0$
and therefore $I$ contains a nonzero constant. The Lemma follows immediately. $\hfill \Box$
\end{pf}

\bigskip
We let 
$\frak{p} = \mathcal{M}\Delta$ be the principal prime ideal in $\mathcal{M}$ generated by $\Delta$.

\begin{lemma}\label{lemma3.2} Suppose that $I \subseteq \mathcal{M}$ is a nonzero
$d$-ideal. Then $I \supseteq \frak{p}^r$ for some integer $r \geq 0$. 
\end{lemma}
\noindent
{\bf Proof:} Let 
\begin{eqnarray*}
J = \{ x \in \mathcal{M} \ | \ \Delta^n x \in I, \ n \geq 0 \}.
\end{eqnarray*}
$J/I$ is  the $\Delta$-torsion submodule of $\mathcal{M}/I$. 
If $\Delta^m \in J$ for some
$m \geq 0$ then the Lemma holds. So we may, and shall, assume that this is 
\emph{not} the case. Now no power $\Delta^m$ is contained in $J$ either. 
Because $J$ is itself a $d$-ideal,  there is no loss in assuming that $I=J$. 
Then $x \in \mathcal{M}, \Delta x \in I \Rightarrow x \in I$, so that (\ref{Dcondition}) holds.
By Lemma \ref{lemma3.1} it follows that $I = \mathcal{M}$, in which case the Lemma is clear.
 \ \ \ $\Box$

\section{Proof of Theorems 1 and 2}
We keep previous notation and assumptions. We separate the main part of the proof of Theorem 1, which is the following.
\begin{thm}\label{thmfingen}
$\mathcal{H}(V)$ is a \emph{finitely generated} $\mathcal{M}$-module.
\end{thm}

First  we show by induction on $\dim V$ that if Theorem \ref{thmfingen} holds for \emph{irreducible} representations then it holds in general. Indeed, let 
\begin{eqnarray*}
0 \rightarrow U \rightarrow V \rightarrow W \rightarrow 0
\end{eqnarray*}
be a short exact sequence of $\Gamma$-modules with $W$ irreducible. By Lemma
\ref{lemmaexact} we have an exact sequence of $\mathcal{M}$-modules
\begin{eqnarray*}
0 \rightarrow \mathcal{H}(U) \rightarrow \mathcal{H}(V) \rightarrow \mathcal{H}(W).
\end{eqnarray*}
By induction $\mathcal{H}(U)$ is finitely generated. Assuming that $\mathcal{H}(W)$ is finitely generated, $\mathcal{H}(V)/\mathcal{H}(U)$ is also finitely generated since it isomorphic to an $\mathcal{M}$-submodule of $\mathcal{H}(W)$. Now the finite generation of $\mathcal{H}(V)$ follows.

\bigskip
For the duration of the proof of Theorem \ref{thmfingen}, we assume that $V$ is an \emph{irreducible} 
$\Gamma$-module of dimension $p$. One consequence of this is that \emph{every} nonzero vector-valued modular form $F(\tau) \in \mathcal{H}(k, V)$
is essential. Fix a nonzero $F(\tau) \in \mathcal{H}(k_0, V)$ (cf. (\ref{graded})), and 
 introduce the graded $\mathcal{M}$-submodule
\begin{eqnarray*}
\mathcal{G} &=& \bigoplus_{k \geq 0} \mathcal{G}(k_0+2k) =  \sum_{i = 0}^{p-1} \mathcal{M}d^iF
\end{eqnarray*}
of $\mathcal{H}(V)$.
 Being linearly independent, the component functions
of $F(\tau)$ cannot satisfy a linear differential equation of order \emph{less} than $p$. This implies that $\mathcal{G}$ is a \emph{direct sum}
\begin{eqnarray}\label{G'dirsum}
\mathcal{G} = \bigoplus_{i=0}^{p-1} \mathcal{M}d^i(F),
\end{eqnarray}
and in particular it is a finitely generated free $\mathcal{M}$-module. The next result is crucial.

\noindent
\begin{lemma} \label{thm7.2} There is a constant $B$, depending only on $V$, such that 
\begin{eqnarray}\label{boundedness}
\dim \left(\mathcal{H}(k_0+2k)/\mathcal{G}(k_0+2k)\right) \leq B
\end{eqnarray}
holds for all $k \geq 0$.
\end{lemma}
\noindent
\begin{pf} Because the action of $d$ on $\mathcal{H}(k, \rho)$ raises weights by $2$, it follows from (\ref{G'dirsum}) that
\begin{eqnarray}\label{G'eq}
\dim \mathcal{G}(k_0+2k) = \sum_{i=0}^{p-1} \dim \mathcal{M}_{2k - 2i}.
\end{eqnarray}
It is well-known that $\dim \mathcal{M}_{2k} = [ k/6]$ or $[k/6]+1$ according as
$2k$ is, or is not,  congruent to $2 \ (\mbox{mod} \ 12)$. Using this, it follows from (\ref{G'eq})
that there is a constant $B'$ depending only on $p$, such that
\begin{eqnarray*}
\left| \dim \mathcal{G}(k_0 + 2k) - pk/6 \right| \leq B'.
\end{eqnarray*}
The Lemma follows from this together with Lemma \ref{lemmaHest}.  $\hfill \Box$
\end{pf}

\bigskip
Consider the tower of $\mathcal{M}$-modules
\begin{eqnarray}\label{Mtower}
\mathcal{H}(V) \supseteq \mathcal{T} \supseteq \mathcal{G}\supseteq 0
\end{eqnarray}
where 
\begin{eqnarray}\label{Tdef}
\mathcal{T} = \{ x \in \mathcal{H}(V) \ | \ \Delta^r x \in \mathcal{G}, \ r \geq 0 \}.
\end{eqnarray}
Thus $\mathcal{T}/\mathcal{G}$ is  the $\Delta$-torsion submodule of $\mathcal{H}(V)/\mathcal{G}$. In order to show that $\mathcal{H}(V)$ is finitely generated, it suffices to show that
the $\mathcal{M}$-modules $\mathcal{H}(V)/\mathcal{T}, \mathcal{T}/\mathcal{G}$, and $\mathcal{G}$ are each finitely generated. The finite generation of $\mathcal{G}$ has already been established.

\begin{lemma} $\mathcal{H}(V)/\mathcal{T}$ is a finitely generated $\mathcal{M}$-module.
\end{lemma}
\begin{pf} By construction, $\Delta$ does not annihilate any nonzero element of the quotient module in question. We assert that $\mathcal{H}(V)/\mathcal{T}$ is a \emph{torsion-free}
$\mathbb{C}[\Delta]$-module. If not, using the fundamental theorem of algebra we can find a nonzero
homogeneous element $x \in \mathcal{H}(V)/\mathcal{T}$ and a scalar $\lambda$ such that $\Delta + \lambda$
annihilates $x$. Because the action of $\Delta$ raises weights by $12$, this forces $\lambda = 0$, so that $\Delta$ annihilates $x$, a contradiction.

\bigskip
Now suppose that the Lemma is false. Then for any integer $n$ we can find a finitely generated 
graded $\mathcal{M}$-submodule of $\mathcal{H}(V)/\mathcal{T}$ generated by no fewer than $n$ elements, and therefore also a finitely generated graded $\mathbb{C}[\Delta]$-submodule 
$\mathcal{I}_n$, say,  generated by no fewer than $n$ elements. As a finitely generated torsion-free
$\mathbb{C}[\Delta]$-module, $\mathcal{I}_n$ is necessarily free because $\mathbb{C}[\Delta]$ is a principal ideal domain. Thus we have shown that 
$\mathcal{H}(V)/\mathcal{T}$ contains finitely generated graded free $\mathbb{C}[\Delta]$-modules of arbitrarily large rank. It is easy to see that this is not consistent with the boundedness of the grading on
$\mathcal{H}(V)/\mathcal{T}$ established in Lemma \ref{thm7.2}, which contradiction completes the proof of the Lemma. $\hfill \Box$
\end{pf}

\bigskip
It remains to prove that $\mathcal{T}/\mathcal{G}$ is a finitely generated $\mathcal{M}$-module.
We prove a bit more than we need, namely
\begin{lemma}\label{thm7.3} There is an integer $r \geq 0$ such that $\frak{p}^r \subseteq$ 
Ann$_{\mathcal{M}}(\mathcal{T}/\mathcal{G})$. Moreover, 
$\mathcal{T}/\mathcal{G}$ is a finitely generated $\mathcal{M}/\frak{p}^r$-module. 
\end{lemma}
\begin{pf}  Define a sequence of graded $\mathcal{M}$-modules $\mathcal{A}_n, \ n = -1, 0, 1, \hdots$ as follows.
$\mathcal{A}_{-1} = \mathcal{G}\Delta$, and for $n \geq 0$, 
\begin{eqnarray*}
\mathcal{A}_n/\mathcal{G} = \mbox{Ann}_{\mathcal{T}/\mathcal{G}}(\frak{p}^n).
\end{eqnarray*}
Then
\begin{eqnarray}\label{Aseq}
\mathcal{A}_{-1}  \subseteq \mathcal{G} = \mathcal{A}_0 \subseteq \mathcal{A}_1 \subseteq \hdots
\end{eqnarray}
is an ascending sequence of $\mathcal{M}$-modules such that each quotient $\mathcal{A}_n/\mathcal{A}_{n-1}$ is annihilated by $\frak{p}$. We assert that $\mathcal{A}_n/\mathcal{A}_{n-1}, n \geq 0,$ is a finitely generated, torsion-free 
$\mathcal{M}/\frak{p}$-module.
If $n=0$ the result holds because $\mathcal{G}$ is a finitely generated,  free $\mathcal{M}$-module. Proceeding by induction on $n$, consider the morphism
\begin{eqnarray*}
&&\varphi: \mathcal{A}_{n+1}/\mathcal{A}_n \rightarrow \mathcal{A}_n/\mathcal{A}_{n-1} \\
&& \ \ \  \ \ \ \  a + \mathcal{A}_n \mapsto \Delta a + \mathcal{A}_{n-1}
\end{eqnarray*}
of $\mathcal{M}/\frak{p}$-modules. By construction, $\varphi$ is an
\emph{injection}.  By induction, $\mathcal{A}_n/\mathcal{A}_{n-1}$ is a finitely generated, torsion-free 
 $\mathcal{M}/\frak{p}$-module, hence the same is true for any submodule, and in particular for
  $\mathcal{A}_{n+1}/\mathcal{A}_n \cong$ im $\varphi$.
 This proves our assertion.
 
 \medskip
 From what we have established so far, it follows that every nonzero quotient
 $\mathcal{A}_{n+1}/\mathcal{A}_n$ contains a graded free $\mathcal{M}/\frak{p}$-submodule. Lemma \ref{thm7.2}
 then implies that the sequence (\ref{Aseq}) stabilizes. Let
 $r\geq 0$ be the least integer such that $\mathcal{A}_r = \mathcal{A}_{r+1}$. Since 
 $\mathcal{T}/\mathcal{G}$ is a $\Delta$-torsion module
it follows that $\mathcal{T} = \mathcal{A}_r$, and this is equivalent to the first assertion of the Lemma. Furthermore, since each $\mathcal{A}_{n+1}/\mathcal{A}_n$
 is finitely generated as $\mathcal{M}/\frak{p}$-module then $\mathcal{T}/\mathcal{G}$ is 
 finitely generated as $\mathcal{M}/\mathfrak{p}^r$-module.  This completes the proof of the
 Lemma, and with it also the proof of  Theorem \ref{thmfingen}. $\hfill \Box$ 
 \end{pf}

 \bigskip
 \bigskip
We now turn to the proof of Theorem \ref{thm1}. In order to show that $\mathcal{H}(V)$ is a
 free $\mathcal{M}$-module it suffices (after Theorem \ref{thmfingen}) to establish that
 $\mathcal{H}(V)$ is  \emph{Cohen-Macaulay}. See, for example, Section 4.3 of
 \cite{B} for further details and  facts that we use below. Write $\mathcal{H} = \mathcal{H}(V)$.
 We assert that $(\Delta, Q)$
is a \emph{regular sequence} for $\mathcal{H}$. This means that the following two conditions hold:
\begin{eqnarray*}
&&(a) \ \mbox{Ann}_{\mathcal{H}}(\Delta) = 0, \ 0 \not=  \Delta \mathcal{H}
\not=  \mathcal{H}, \\
&&(b) \ \mbox{Ann}_{\mathcal{H}/\Delta \mathcal{H}}(Q) = 0,   \ 
0 \not= Q (\mathcal{H}/\Delta \mathcal{H}) \not= \mathcal{H}/\Delta \mathcal{H}. 
\end{eqnarray*}
These facts follow easily from a consideration of the $q$-expansions of components of vector-valued modular forms in $\mathcal{H}$. In more detail,  (a) holds because $\mathcal{H}$ is  torsion-free
as $\mathcal{M}$-module and because multiplication by $\Delta$ raises weights by $12$. As for
(b), notice that $\Delta \mathcal{H}$ consists of those holomorphic vector-valued modular forms
$F(\tau) = (f_1(\tau), \hdots, f_p(\tau))^t$ such that the coefficient of the leading power $q^{m_j}$
of the $q$-expansion of $f_j(\tau)$ 
(cf. (\ref{compqexp})) \emph{vanishes}, $1 \leq j \leq p$. On the other hand, multiplication by $Q$ raises weights by $4$ and does not change the order of vanishing at $\infty$. (b) follows immediately.

\bigskip
Because $\mathcal{M}$ has Krull dimension $2$, the existence of the regular sequence 
$(\Delta, Q)$
of length $2$ means that $\mathcal{H}$ is indeed Cohen-Macaulay. Then because 
$\mathcal{M}=\mathbb{C}[Q, R]$ is a weighted polynomial algebra, $\mathcal{H}$ is a finitely generated
free $\mathcal{M}$-module.

\bigskip
Finally, we have to show that the rank of $\mathcal{H}$ as $\mathcal{M}$-module is exactly $p$.
The proof is similar to that of Lemma \ref{thm7.2}. Namely, suppose that
$F_1, \hdots, F_l$ are vector-valued modular forms of weights $k_0+2k_1, \hdots, k_0+2k_l$
respectively and also free generators of $\mathcal{H}$ as $\mathcal{M}$-module. Then 
\begin{eqnarray*}
\dim \mathcal{H}(k_0+2k, \rho) = \sum_{j=1}^l \dim \mathcal{M}_{2k-2k_j},
\end{eqnarray*}
so that for $k \gg 0$ we have
\begin{eqnarray*}
 \sum_{j=1}^l [(k-k_j)/6] \leq \dim \mathcal{H}(k_0+2k, \rho)  \leq \sum_{j=1}^l [(k-k_j)/6] + 1.
\end{eqnarray*}
Thus there is a constant $C$ such that for all $k$ we have
\begin{eqnarray*}
|\dim \mathcal{H}(k_0+2k, \rho) - lk/6| < C.
\end{eqnarray*}
Comparing this with Lemma \ref{lemmaHest} shows that $l=p$. This completes the proof of
Theorem \ref{thm1}.

\bigskip
We turn to the proof of Theorem 2. Part (a) is included in Lemma \ref{lemmaexact}.
We restate part (b) as follows:

 \begin{lemma}\label{lemmasplitinj} Suppose that $i: U \rightarrow V$ is an injective 
 morphism of $\Gamma$-modules. Then 
 $\mathcal{H}(U) \stackrel{i_*}{\rightarrow} \mathcal{H}(V)$
is a \emph{split injection}. That is, $i_*\mathcal{H}(U)$ is a direct summand of $\mathcal{H}(V)$
as $\mathcal{M}$-module.
 \end{lemma}
 \begin{pf} We know from Lemma \ref{lemmaexact} that $i_*$ is an injection. 
  In the following, we identify $\mathcal{H}(U)$ with its $i_*$-image in $\mathcal{H}(V)$ and $U$ with its $i$-image in $V$. First we prove that $\mathcal{H}(U)$ is a \emph{pure} $\mathcal{M}$-submodule of $\mathcal{H}(V)$, ie., $\mathcal{H}(V)/\mathcal{H}(U)$
  is \emph{torsion-free}. 
  Suppose that $0 \not= g(\tau) \in \mathcal{M}_{2k}, F(\tau) \in \mathcal{H}(k_1, V)$ and
  $gF \in \mathcal{H}(U)$. We have to show that $F \in \mathcal{H}(U)$.
  If $\tau \in \mathcal{H}$ we have 
\begin{eqnarray*}
g(\tau)F(\tau) \in U \subseteq V.
\end{eqnarray*}
If $g(\tau) \not= 0$ this implies that $F(\tau) \in U$. Since $g$ is nonzero, the zeros of $g(\tau)$
in $\mathcal{H}$ are \emph{discrete}. Since $F$ is continuous and $U$ closed in $V$, it follows that
$F(\mathcal{H}) \subseteq U$, that is $F \in \mathcal{H}(U)$. This establishes the purity of
$\mathcal{H}(U)$. 

\medskip
Let $\mathcal{M}^+$ be the maximal ideal of $\mathcal{M}$ generated by $Q$ and $R$. From Theorem \ref{thm1} we know that 
\begin{eqnarray*}
E = \mathcal{H}(V)/\mathcal{M}^+\mathcal{H}(V)
\end{eqnarray*}
is a finite-dimensional $\mathbb{C}$-linear space of dimension $p$, and that a
set $F_1, \hdots, F_p$ of $p$ (homogeneous) elements of $\mathcal{H}(V)$ is a set of free generators
(as $\mathcal{M}$-module) if, and only if, $F_1, \hdots, F_p$ maps onto a 
\emph{basis} of $E$ via the canonical projection $\mathcal{H}(V) \rightarrow \mathcal{H}(V)/\mathcal{M}^+\mathcal{H}(V).$
Because of the purity of $\mathcal{H}(U)$ it is easy to see that 
\begin{eqnarray*}
\mathcal{M}^+\mathcal{H}(V) \cap \mathcal{H}(U) = \mathcal{M}^+\mathcal{H}(U).
\end{eqnarray*}
Therefore, a second application of Theorem \ref{thm1} shows that
\begin{eqnarray*}
E_0 = (\mathcal{H}(U) + \mathcal{M}^+\mathcal{H}(V))/\mathcal{M}^+\mathcal{H}(V)
\cong \mathcal{H}(U)/\mathcal{M}^+\mathcal{H}(U) 
\end{eqnarray*}
is a $\mathbb{C}$-linear subspace of $E$ of dimension $n = \dim U$. Choose
homogeneous elements $G_1, \hdots, G_n \in \mathcal{H}(U)$ which map onto a basis of
$E_0$, and extend this set to a basis of $E$ by choosing appropriate homogeneous elements
$G_{n+1}, \hdots, G_p$ in $\mathcal{H}(V)$. It follows that 
\begin{eqnarray*}
\mathcal{H}(V) = \mathcal{H}(U) \oplus N
\end{eqnarray*}
where $N$ is the $\mathcal{M}$-submodule generated by $G_{n+1}, \hdots, G_p$. This completes the proof of the Lemma.
$\hfill \Box$
 \end{pf}
 
\bigskip
Turning to the proof of Theorem 2(c), set $\mathcal{H}' = j_*\mathcal{H}(V)$. After part (b) we know that there is a direct sum decomposition of $\mathcal{M}$-modules
\begin{eqnarray*}
\mathcal{H}(V) = i_*\mathcal{H}(U) \oplus \mathcal{H}',
\end{eqnarray*}
and by Theorem 1, $\mathcal{H}(V)$ and $\mathcal{H}(U)$ are free  of ranks $\dim V$ and $\dim U$ respectively. Then $\mathcal{H}'$ is free of rank $\dim V$-$\dim U$.

\medskip
It remains to show that $\mathcal{H}(W)/\mathcal{H}'$, which is a graded $\mathcal{M}$-module,  has finite dimension. In other words, we must show that for all large enough $k$,
$\mathcal{H}(k, W) \subseteq \mathcal{H}'$.
Let us assume, as we may, that $V$ furnishes a representation $\rho$ of $\Gamma$ which is 
\emph{upper triangular}. More precisely, for $\gamma \in \Gamma$ set
\begin{eqnarray}\label{ut}
\rho(\gamma) = 
\left(\begin{array}{cc} \alpha(\gamma) & \beta(\gamma) \\0 & \sigma(\gamma)\end{array}\right).
\end{eqnarray}
Then $\alpha$ and $\sigma$ are representations of $\Gamma$ corresponding to the $\Gamma$-modules
$U$ and $V/U \cong W$ respectively. Suppose that $r= \dim U, s = \dim W$, and that
\begin{eqnarray*}
F(\tau) = (f_1(\tau), \hdots, f_r(\tau), f_{r+1}(\tau), \hdots, f_{r+s}(\tau))^t
\end{eqnarray*}
is an element in $\mathcal{H}(k, V)$ adapted to the choice of basis for which $\rho$
is as in (\ref{ut}). Then $F'(\tau) = (f_r(\tau), f_{r+1}(\tau), \hdots, f_{r+s}(\tau))^t$ is an element of
$\mathcal{H}(k, W)$ and the morphism $j_*: \mathcal{H}(V) \rightarrow \mathcal{H}(W)$ 
induced by $j: V \rightarrow W$ is just the map
\begin{eqnarray}\label{jmorph}
F(\tau) \mapsto F'(\tau).
\end{eqnarray}

\medskip
Now the operator $D$ on $\mathcal{H}(V)$ (and $\mathcal{H}(W)$) acts in a componentwise
fashion. Then the previous discussion makes it clear that $j_*$
\emph{is a morphism of graded $\mathcal{R}$-modules}. As a result, 
$\mathcal{H}'$ and $\mathcal{H}(W)/\mathcal{H}'$ are both graded $\mathcal{R}$-modules.
By Theorem 1, $\mathcal{H}(W)$ is a free $\mathcal{M}$-module of rank $s$. 
Another application of Theorem 1, together with Theorem 2(b), shows that
$\mathcal{H}'$ is also a free $\mathcal{M}$-module of rank $s$. It follows (cf. the proof of
Lemma \ref{thm7.2}) that the grading on $\mathcal{H}(W)/\mathcal{H}'$ has
\emph{bounded dimension} in the sense that there is an upper bound on the dimension of the
homogeneous subspaces. 

\medskip
Consider $\mathcal{H}(W)/\mathcal{H}'$ as $\mathcal{M}$-module. Because
it is finitely generated, the boundedness of the dimensions of the homogeneous subspaces implies that
the annihilator Ann$_{\mathcal{M}}(\mathcal{H}(W)/\mathcal{H}')$ is a
\emph{nonzero} ideal, call it $J$. Moreover, because $\mathcal{H}(W)/\mathcal{H}'$
is a graded $\mathcal{R}$-module then $J$ is a graded $d$-ideal in the sense of
Lemma \ref{lemma3.1}. By Lemma \ref{lemma3.2} we conclude that there is an integer $t$ such that
\begin{eqnarray}\label{Dtors}
J \supseteq \frak{p}^t.
\end{eqnarray}

We will prove the following assertion: 
\begin{eqnarray}\label{noD}
&&\mbox{for all large enough $k$, the homogeneous
subspace}   \notag \\
&&\mbox{$(\mathcal{H}(W)/\mathcal{H}')_k$ of weight $k$ contains
\emph{no nonzero $\Delta$-torsion}.}
\end{eqnarray}
If this is so, it consistent with (\ref{Dtors}) only if $(\mathcal{H}(W)/\mathcal{H}')_k = 0$
for large enough $k$. This is equivalent to the containment $\mathcal{H}(k, W)
\subseteq \mathcal{H}'$, which is what we seek to prove. Therefore, the proof of
Theorem 2(c) is reduced to establishing (\ref{noD}).

\medskip
Suppose that $G(\tau) \in \mathcal{H}(k, W)$ satisfies $\Delta G \in \mathcal{H}'$.
From our discussion of (\ref{jmorph}), this means that there is
$F \in \mathcal{H}(k+12, V)$ with
\begin{eqnarray*}
F(\tau) = (f_1(\tau), \hdots, f_r(\tau), \Delta(\tau)G(\tau))^t.
\end{eqnarray*}
We now apply Lemma \ref{lemmamus}. It tells us that for large enough $k$, we can always find
$H(\tau)=(h_1(\tau), \hdots, h_r(\tau))^t  \in \mathcal{H}(k+12, U)$ such that the leading term
of $h_j(\tau)$ \emph{coincides} with that of $f_j(\tau)$ for $1 \leq j \leq r$. Now consider
the vector-valued modular form
\begin{eqnarray*}
F^0(\tau) = F(\tau) - (H(\tau), \underbrace{0, \hdots, 0}_s)^t.
\end{eqnarray*}
By construction we have $\Delta^{-1}F^0(\tau) \in \mathcal{H}(k, V)$. Moreover,
\begin{eqnarray*}
j_*(\Delta^{-1}F^0(\tau)) = \Delta^{-1} j_*(F(\tau)) = G(\tau),
\end{eqnarray*}
so that in fact $G(\tau) \in \mathcal{H}'$. This proves that (\ref{noD}) holds
for all large enough $k$, and the proof of Theorem 2(c) is complete. $\hfill \Box$

\section{Proof of Theorems 3 and 4}
For the first two results, we fix an \emph{essential} vector-valued modular form 
$F \in \mathcal{H}(k_1, V)$. 
 In this context we continue to use the notation of the previous Section. In particular, we set
 \begin{eqnarray*}
\mathcal{G} = \bigoplus_{i=0}^{p-1} \mathcal{M}d^iF.
\end{eqnarray*}
Because $F$ is essential the last display is indeed, as in (\ref{G'dirsum}), a direct sum of $\mathcal{M}$-modules.
The argument is the same  as before. We also have the corresponding tower of $\mathcal{M}$-modules (\ref{Mtower}) with $\mathcal{T}$ as in (\ref{Tdef}). Introduce
\begin{eqnarray*}
\mathcal{G}' = \sum_{i = 0}^{\infty} \mathcal{M}d^iF.
\end{eqnarray*}
$\mathcal{G}'$ is the cyclic $\mathcal{R}$-submodule of $\mathcal{H}(V)$ generated by $F$. Of course it contains
$\mathcal{G}$, and generally the containment is \emph{proper}.

\begin{lemma}\label{thmH/G'} The $\mathcal{M}$-module $\mathcal{H}(V)/\mathcal{G}'$
is a $\Delta$-torsion module, i.e., it is annihilated by some power of $\Delta$.
\end{lemma}
\begin{pf} The proof is similar to a part of the argument used in the proof of
Theorem 2(c). Briefly, it goes as follows.
 Arguing as in Lemma \ref{thm7.2} we find that the homogeneous spaces
of  $\mathcal{H}/\mathcal{G}'$ are of bounded dimension. Because of finite generation, it follows that
the annihilator $I = Ann_{\mathcal{M}}(\mathcal{H}/\mathcal{G}')$ is nonzero. Because
 $\mathcal{H}/\mathcal{G}'$ is an $\mathcal{R}$-module then $I$ is a nonzero $d$-ideal of
 $\mathcal{M}$. By Lemma \ref{lemma3.2} it follows that $\frak{p}^r \subseteq I$ for some $r \geq 1$,
 as required.  $\hfill \Box$
\end{pf}

\begin{lemma}\label{thm9.1} Suppose that $F = (f_1, \hdots, f_p)^t$ and that the component
functions
$f_j$ have $q$-expansions 
\begin{eqnarray*}
f_j(\tau) = c_j q^{\lambda_j} + \hdots, \ c_j \not= 0, \ 1 \leq j \leq p.
\end{eqnarray*}
Then the following are equivalent:
\begin{eqnarray*}
 &&(a) \  \mathcal{T} = \mathcal{G},  \\
 &&(b) \ m_j = \lambda_j \  \mbox{for $1 \leq j \leq p$ and no two of the $m_j$ are equal}.
 \end{eqnarray*}
 (Here, the $m_j$ are as in (\ref{Tdiag}).)
\end{lemma}
\begin{pf} We keep the notation used in the preceding proof. 
 We may, and shall, assume that the component functions $f_1(\tau), \hdots, f_p(\tau)$
have normalized $q$-expansions
\begin{eqnarray*}
f_j(\tau) = q^{\lambda_j} + \hdots, \ j = 1, \hdots, p
\end{eqnarray*}
(recall that $m_j\leq \lambda_j$ and $\lambda_j - m_j \in \mathbb{Z}$), with
\begin{eqnarray}\label{aorder}
&& \ 0 \leq \lambda_1 \leq \hdots \leq \lambda_s < 1 \leq \lambda_{s+1} \leq \hdots \leq \lambda_p. 
\end{eqnarray}

Let $F' \in \mathcal{G}$ be a vector-valued modular form of weight $k_1+2k$, with
\begin{eqnarray*}
 F' &=& \sum_{i=0}^{p-1} g_id^iF, \ \ g_i \in \mathcal{M}_{2k - 2i}, \\
g_i &=& a_i + O(q), \ a_i \in \mathbb{C}, \ 0 \leq i \leq p-1.
\end{eqnarray*}
Because the $d^iF, \ 0 \leq i \leq p-1,$ are linearly independent over $\mathcal{M}$,
it follows that 
$\mathcal{T}  \not= \mathcal{G}$ if, and only if, it is possible to choose $a_i$'s not
all equal to $0$ in such a way that the component functions of $F'$ nevertheless all vanish to order at least
$1$ at $\infty$.

\medskip
The vector consisting of leading coefficients of the components of $F'$ is equal to $Av^t$ where
$v = (a_0, \hdots, a_{p-1})$ and $A$ is the $p \times p$ matrix
\begin{eqnarray*}
A = \left(\begin{array}{cccc}1 & \lambda_1-1/12 & \hdots  &  \\1 & \lambda_2 - 1/12 & \hdots  &  \\ \vdots & \vdots &   & \vdots \\1 & \lambda_p - 1/12 & \hdots  &  \end{array}\right).
\end{eqnarray*}
We easily see that $A$ is similar to the Vandermonde matrix
\begin{eqnarray*}
 \left(\begin{array}{cccc}1 & \lambda_1 & \hdots & \lambda_1^{p-1} \\1 & \lambda_2 & \hdots & 
\lambda_2^{p-1} \\ \vdots & \vdots &  & \vdots \\ 1 & \lambda_p & \hdots & \lambda_p^{p-1}  \end{array}\right), 
\end{eqnarray*}
and in particular, $A$ is invertible if, and only if, all of the $\lambda_i$'s
are \emph{distinct}.

\medskip
Write $A$ in block form 
\begin{eqnarray*}
A = \left(\begin{array}{cc}  U  & V\\ W & X   \end{array}\right)
\end{eqnarray*}
where $U$ is an $s \times s$ matrix and $s$ is as in (\ref{aorder}). In order to be able to choose the vector
$v$ of $a_i$'s so that the components of $\sum g_i d^iF$ vanish to order at least $1$ at $\infty$, it is necessary and sufficient
to solve the system of equations
\begin{eqnarray}\label{ve}
\left(\begin{array}{cc} U&V\\ W&X   \end{array}\right) \left( \begin{array}{c} v_1 \\ v_2
\end{array} \right) =  \left( \begin{array}{c} 0 \\ * \end{array} \right) 
\end{eqnarray}
Here, $v^t = (v_1, v_2)^t$ and $*$ is arbitrary. This is because the `lower piece' of $v^t$ corresponds to those
$\lambda_j \geq 1$ and no condition is imposed on the corresponding coefficients.

\medskip
Now suppose that (b) of Theorem \ref{thm9.1} holds. Then 
$A$ is invertible, $s=p$, and the only way to solve (\ref{ve}) is with
$v = 0$. So $\mathcal{T} = \mathcal{G}$ in this case and (a) holds.

\medskip
Suppose that the $\lambda_j$'s are 
\emph{not} distinct.  Then $A$ is singular, and we may take $v$ to be any nonzero vector annihilated
by $A$. So in this case (a) does not hold. On the other hand, if the $\lambda_j$'s
are distinct but not all less than $1$ (i.e. $s<p$),  then we may choose $*$ in (\ref{ve}) to be nonzero and
take $v^t = A^{-1}(0, *)^t \not = 0.$ So (a) does not hold in this case either. 
This shows that if (b) is false then so is (a), and the proof of the Lemma is complete. $\hfill \Box$
\end{pf}

 \bigskip
We now prove Theorem 3. Suppose first that $\mathcal{H}(V) = \mathcal{R}F$
is a cyclic $\mathcal{R}$-module with generator $F$. In our earlier notation, we have
$\mathcal{H}(V) = \mathcal{G}'$. 
Let the components of $F$ be as in Lemma \ref{thm9.1}. Note that $F$ is necessarily an
essential vector-valued modular form. Let the weight of $F$ be $k_0$.
 
 \medskip
 We claim that $\mathcal{G} = \mathcal{G}'$. We argue as follows.
 By Theorem \ref{thm1}, $\mathcal{H}(V)$ has an $\mathcal{M}$-basis of cardinality $p$. 
 Because $\mathcal{H}(V) = \mathcal{R}F$, up to scalars $F$ is the unique
 nonzero vector-valued modular form in $\mathcal{H}(V)$ of least weight. Thus, we may 
 choose $F$ to be a member of a basis of $\mathcal{H}(V)$. Let $r$ be the maximal integer
 such that $F, dF, \hdots , d^rF$ is part of a basis, with
 \begin{eqnarray*}
N = \sum_{i=0}^r \mathcal{M}d^iF
\end{eqnarray*}
the $\mathcal{M}$-submodule of $\mathcal{H}(V)$  spanned by these elements.
 Thus $0 \leq r \leq p-1$, and we want to show that $r=p-1$.  If $d^{r+1}F \in N$ then there is a linear relation
 \begin{eqnarray*}
d^{r+1}F = \sum_{j=1}^r g_jd^jF, \ \ g_j \in \mathcal{M}.
\end{eqnarray*}
Since $F$ is essential, such a relation implies that $r+1\geq p$ and we are done. If
$d^{r+1}F \notin N$ then $d^{r+1}F +N$ is,  up to scalars, the unique nonzero element 
 of least weight in the free $\mathcal{M}$-module $\mathcal{H}(V)/N$. 
 Then (up to scalars) $d^{r+1}F+N$ is necessarily a member of any
 $\mathcal{M}$-basis of $\mathcal{H}(V)/N$, in which case $F, dF, \hdots, d^{r+1}F$
is part of a basis of $\mathcal{H}(V)$. This contradicts the definition of $r$, and proves
our claim.

\medskip
Having established that $\mathcal{G} = \mathcal{G}'$, it follows that
\begin{eqnarray*}
\mathcal{H}(V) = \mathcal{G}.
\end{eqnarray*}
In particular, we have $\mathcal{T} = \mathcal{G}$, so that
Lemma \ref{thm9.1} is applicable and we can conclude that (b) of
Lemma \ref{thm9.1} holds. Furthermore, since $d^pF \in \mathcal{G}$ there is a relation
\begin{eqnarray*}
d^pF + \sum_{j=2}^{p-1} g_jd^jF, \ g_j \in \mathcal{M}_{2p-2j}.
\end{eqnarray*}
This defines a MLDE $L_{k_0}[f] = 0$ where $L$ is the differential operator
\begin{eqnarray*}
L = D^p + \sum_{j=2}^{p-1} g_jD^j. 
\end{eqnarray*}
The roots of the corresponding indicial equation are the $\lambda_j$. As we have seen,
these coincide with the $m_j$, and are distinct. That the weight $k_0$ is determined by
(\ref{wtrestr}) is proved in \cite{M1}. We have now established all of the conditions stated in
Theorem 3 under the assumption that $\mathcal{H}(V) = \mathcal{R}F$.
 
 \bigskip
 As for the converse, suppose that we have a MLDE satisfying the conditions of Theorem 3.
 The solution space $V$ defines an element $F \in \mathcal{H}(k_0, V)$
 (cf. \cite{M1}), and since the roots of the indicial equation are real and distinct
 then the monodromy matrix $\rho(T)$ is unitarizable. (Here, $V$ affords the representation
 $\rho$ of $\Gamma$.)
 
 \medskip
 Let $\mathcal{G}, \mathcal{G}'$ have the same meaning as before, where now
 $F$ is determined by the MLDE as in the previous paragraph. Thus in fact $\mathcal{G}= \mathcal{G}'$, and we 
 have to show that $\mathcal{H}(V) = \mathcal{G}'$.
 To see this, note that in the present situation part (b) of Lemma \ref{thm9.1}
 holds. That result then shows that $\mathcal{T} = \mathcal{G}$, i.e.,
 $\mathcal{H}(V)/\mathcal{G}$ contains no nonzero $\Delta$-torsion. On the other hand,
 $\mathcal{H}(V)/\mathcal{G}'$ is a $\Delta$-torsion module by Lemma \ref{thmH/G'}.
 Since $\mathcal{G} = \mathcal{G}'$, the only way to reconcile these statements is the conclusion that
 $\mathcal{H}(V) = \mathcal{G}$. This completes the proof of Theorem 3. $\hfill \Box$

\bigskip
Finally, we consider Theorem 4.
If $V$ is not indecomposable it is a direct sum of two $1$-dimensional $\Gamma$-modules,
so we start with the $1$-dimensional case. 
Of course this  is well-known, but it is interesting
nonetheless to reconsider it from our current perspective.
In this case $\mathcal{H}(V)$ is a free $\mathcal{M}$-module of rank $1$. Thus if $F_0$
is a nonzero vector-valued modular form of minimal weight $k_0$ then
\begin{eqnarray*}
\mathcal{H}(V) =\mathcal{M}F_0 = \mathbb{C}F_0 \oplus \mathbb{C}QF_0 \oplus \hdots
\end{eqnarray*}
Since there are no nonzero vector-valued forms of weight $k_0+2$ then $DF_0 = 0$.
Since $F_0 = (f_0)$ where $f_0(\tau)$ is a classical modular form of weight $k_0$, the
condition $Df_0 = 0$ implies that $f_0$ is a scalar multiple of $\eta(\tau)^{2k_0} = q^{k_0/12} + \hdots$.
If $k_0 \geq 12$ then $\Delta^{-1}F_0$ is a nonzero vector-valued form of weight \emph{less}
than $k_0$, a contradiction. Thus $0 \leq k_0 \leq 11$, and we arrive at the $12$ possibilities
\begin{eqnarray*}
\mathcal{H}(V) = \mathcal{M} \eta(\tau)^{2k_0}, \ \ 0 \leq k_0 \leq 11
\end{eqnarray*}
corresponding to the various $1$-dimensional representations $\chi^{k_0}$ of $\Gamma$.

\bigskip
We turn to the case that $V$ is indecomposable but not irreduicble.
It will be useful to record some facts about these modules, which are more-or-less proved in \cite{M2}, Theorem 3.1.

\begin{lemma}\label{2dind} Suppose that $V$ is a $2$-dimensional indecomposable module which is not irreducible and furnishes a representation $\rho$ of $\Gamma$. 
Then there is an ordered pair of $12$th roots of unity $(\mu_1, \mu_2)$ such that
$\mu_1/\mu_2$ is a \emph{primitive} sixth root of unity and such that $\rho(T)$ is similar to the diagonal matrix
\begin{eqnarray*}
\left(\begin{array}{cc} \mu_1 & 0 \\0 & \mu_2\end{array}\right).
\end{eqnarray*}
The representation $\rho$ is \emph{characterized} up to equivalence by $(\mu_1, \mu_2)$. Thus there are just $24$ equivalence classes of such representations.  
$\hfill \Box$
\end{lemma}

\medskip
In what follows we take $V$ with $\rho(T)$ diagonal as in the last Lemma. We may, and shall, assume
that $\rho$ is \emph{upper triangular}.
  Let $F_0 \in \mathcal{H}(V)$ be a nonzero vector-valued modular form of least weight $k_0$, say,
and set $F_0 = (f_1(\tau), f_2(\tau))^t$. We first consider the case
\begin{eqnarray}\label{DFnot0}
DF_0 = 0,
\end{eqnarray}
and assume this until further notice. Then $Df_1 = Df_2 = 0$, and since
$f_1, f_2$ are both solutions of the same differential equation they can differ only by an overall scalar.
If they are both \emph{nonzero} then they have $q$-expansions $f_j(\tau) = c_jq^t + \hdots$
with nonzero scalars $c_j$, and from this it follows that $\rho(T)$ is a \emph{scalar matrix}. Because
$\rho$ is indecomposable this is not possible, and we conclude that one of the components of
$F_0$ vanishes identically. Because $\rho$ is upper triangular then $f_2 = 0$. Now $f_1(\tau)$ is a classical modular form of weight $k_0$ which is annihilated by $D$. Thus $f_1$ is a scalar multiple of
$\eta(\tau)^{2k_0}$, and we may take
\begin{eqnarray*}
F_0 = (\eta(\tau)^{2k_0}, 0)^t. 
\end{eqnarray*}
Because $F_0$ has minimal weight then $\Delta^{-1}F_0$ cannot be holomorphic.
It follows that
$0 \leq k_0 < 12$. 

\medskip
We know that $\mathcal{H}(V)$ is free of rank $2$ as $\mathcal{M}$-module, and that
we may take $F_0$ as one of the free generators. Let $k_1$ be the weight of the second free generator
$G = (g_1, g_2)^t$, say. Thus $k_0 \leq k_1$. 
Because $DG$ has weight $k_1+2$ and is not a free generator of $\mathcal{H}(V)$ then
\begin{eqnarray*}
DG = \alpha F_0 = (\alpha \eta^{2k_0}, 0)^t
\end{eqnarray*}
 for some $\alpha \in \mathcal{M}_{k_1 - k_0+2}$. In particular, $Dg_2 = 0$. Because 
 $\rho$ is upper triangular then $g_2$ is a classical holomorphic modular form, and since
 it is also annihilated by $D$ then 
 \begin{eqnarray*}
 g_2 = u\eta(\tau)^{2k_1}
 \end{eqnarray*}
  for a scalar $u$.
We claim that $u \not= 0$. Otherwise, because $\mathcal{H}(V) = \mathcal{M}F_0 + \mathcal{M}G$ it follows that the second component of \emph{every} element in $\mathcal{H}(V)$ vanishes, and in particular there is no
essential vector-valued modular form. This contradicts Lemma \ref{essvvform}. So we may choose
$u=1$. 

\medskip
We next assert that $k_1 \geq 1$. If not, $k_0 = k_1 = 0$ and $DG = 0$. But then
as before  the component functions of $G$ are linearly dependent, and we easily see in this case
that $\rho$ is the direct sum of a pair of $1$-dimensional representations, contradiction.

\medskip
We claim also that $k_1 \leq 11$. Otherwise, consider the vector-valued modular form
\begin{eqnarray*}
H(\tau) &=&   x E_{k_1-k_0+4}(\tau)F_0  + Q G \\
&=& x E_{k_1-k_0+4}(\tau) \left(\begin{array}{c} \eta(\tau)^{2k_0} \\0\end{array}\right)  + Q \left(\begin{array}{c}g_1 \\ \eta(\tau)^{2k_1}\end{array}\right) \in \mathcal{H}(\rho, k_1+4),
\end{eqnarray*}
$x \in \mathbb{C}$. We can choose $x$ so that both components of $H(\tau)$ vanish to order
at least $1$ at $\infty$, so that $0 \not= \Delta^{-1}H(\tau) \in \mathcal{H}(\rho, k_1-8)$.
This forces $\Delta^{-1}H(\tau) \in \mathcal{M}F_0$, in which case
$H(\tau) \in \mathcal{M}F_0$, contradiction.  

\medskip
Since $f_1(\tau) = \eta(\tau)^{2k_0}$ and $g_2(\tau) = \eta(\tau)^{2k_1}$ with
$0 \leq k_0 \leq k_1 \leq 11$ then $m_1 = k_0/12, m_2 = k_1/12$ (notation as in (\ref{Tdiag})).
By Lemma \ref{2dind} we find that 
\begin{eqnarray}\label{m2<m1}
m_2 - m_1 = 1/6 \ \mbox{or} \ 5/6.
\end{eqnarray}
There are exactly $12$ choices of pairs $(m_1, m_2)$ satisfying (\ref{m2<m1}) together with
$0 \leq m_1 , m_2 < 1$. They correspond to half of the indecomposable $\Gamma$-modules described
in Lemma \ref{2dind}.

\bigskip
Now we assume that (\ref{DFnot0}) does \emph{not} hold. Then (up to scalars), the unique pair of homegeneous
free generators for $\mathcal{H}(V)$ as $\mathcal{M}$-module consists
of $F_0$ and $DF_0$. Thus in this case $\mathcal{H}(V) = \mathcal{R}F_0$, and Theorem
3 applies. Thus we know that we may take

\begin{eqnarray*}
F_0(\tau) = \left(\begin{array}{c}f_1(\tau)  \\ f_2(\tau) \end{array}\right)
= \left(\begin{array}{c} q^{m_1} + \hdots \\ q^{m_2} + \hdots\end{array}\right)
\end{eqnarray*}
with
\begin{eqnarray*}
k_0 = 6(m_1+ m_2)-1.
\end{eqnarray*}

\medskip
Because 
$\rho$ is \emph{upper triangular}, $f_2(\tau)$ is a nonzero classical holomorphic modular form
of weight $k_0$ and
$f_2(\tau) = \eta(\tau)^{24m_2} \alpha$
for some $\alpha \in \mathcal{M}_{k_0 - 12m_2}$. Then
\begin{eqnarray*}
0 \leq k_0-12m_2 = 6(m_1- m_2)-1<5,
\end{eqnarray*}
and since $k_0-12m_2$ is the weight of $\alpha$ then it must be $0$ or $4$. Therefore,
\begin{eqnarray}\label{m1<m2}
m_1 - m_2 =  1/6 \ \mbox{or} \ 5/6.
\end{eqnarray}
There are exactly $12$ ordered pairs $(m_1, m_2)$ satisfying (\ref{m1<m2})
together with $0 \leq m_2 < m_1 < 1$. Notice that these correspond, as they must, to
the $12$ classes of indecomposables  \emph{not} covered by (\ref{m2<m1}).

\bigskip
Let us consider the Hilbert-Poincar\'{e} series of $\mathcal{H}(V)$.
By Theorem 2 there is a containment of $\mathbb{Z}$-graded $\mathcal{M}$-modules
\begin{eqnarray}\label{Hcontain}
\mathcal{H}(V) \subseteq \mathcal{H}(\chi^{a}) \oplus \mathcal{H}(\chi^{b}).
\end{eqnarray}
 Here, notation is as in the statement of Theorem 4.

\medskip
 Consider first the case that (\ref{DFnot0}) holds. Here we showed above that
 the two fundamental weights are $k_0=12m_1$ and $k_1 = 12m_2$.
  It follows that, in the notation of (\ref{Hcontain}),
 we have $a = k_0$ and $b= k_1$. Then both $\mathcal{H}(V)$ and
 $\mathcal{H}(\chi^a)\oplus \mathcal{H}(\chi^b)$ have the \emph{same} Hilbert-Poincar\'{e}
 series
 \begin{eqnarray}\label{abPS}
\frac{t^a+t^b}{(1-t^4)(1-t^6)},
\end{eqnarray}
and in particular (\ref{Hcontain}) is an \emph{equality} in this case.

 \medskip
 Now assume  that (\ref{DFnot0}) does \emph{not} hold. 
 We saw that an $\mathcal{M}$-basis of $\mathcal{H}(V)$ has the form $F_0, DF_0$, so that 
\begin{eqnarray*}
PS\ \mathcal{H}(V) =  \frac{t^{k_0}(1+t^2)}{(1-t^4)(1-t^6)} 
\end{eqnarray*}
with $k_0 = 6(m_1+m_2)-1$. Suppose first that $k_0 = 12m_2$, so that $m_1-m_2=1/6$.
In this case, $f_2 = \eta(\tau)^{24m_2}$ has weight $k_0$ and corresponds to the character
$\chi^b$. Thus $b=12m_2=k_0$ and $a = 12m_1 = 12m_2+2 = k_0+2$. So in this case, the Poincar\'{e}
series of  $\mathcal{H}(V)$ and  $\mathcal{H}(\chi^a)\oplus \mathcal{H}(\chi^b)$
again coincide with (\ref{abPS}) and (\ref{Hcontain}) is an equality.

\medskip
 In the remaining case,
when $k_0 = 12m_2+4$ and $m_1-m_2=5/6$, we have
\begin{eqnarray*}
&& PS(\mathcal{H}(\chi^a) \oplus \mathcal{H}(\chi^b)) - PS \ \mathcal{H}(V) \\
=&& \frac{t^{12m_1}+t^{12m_2}}{(1-t^4)(1-t^6)} -  \frac{t^{k_0}(1+t^2)}{(1-t^4)(1-t^6)} \\
=&&  \frac{t^{12m_2}(1+t^{10}-t^4-t^6)}{(1-t^4)(1-t^6)} \\
=&& t^{12m_2}.
\end{eqnarray*}
Thus in this case, the containment (\ref{Hcontain}) is proper and the codimension exactly $1$,
occuring in weight $b=12m_2$. The statement of Theorem 4 in the Introduction is just a reorganization of these calculations.

 \bigskip
 \noindent
 Author's address\\
 Department of Mathematics,  University of California at Santa Cruz
CA 95064 \\
cem@ucsc.edu \\
gem@cats.ucsc.edu

\end{document}